\documentclass[12pt]{article}

\renewcommand{\Re}{\mathrm{Re}}
\newcommand{\C}{\mathbb{C}}  
\renewcommand{\P}{{\bf P}}
\newcommand{\N}{\mathcal{N}}
\newcommand{\interior}{\mathrm{int}}

\usepackage{amsfonts}
\usepackage{amsmath,amsthm}
\usepackage{graphicx}

\title{Maximizing the Closed Loop Asymptotic Decay Rate for the
  Two-Mass-Spring Control Problem}

\author{Didier Henrion$^{1,2}$ \and Michael L. Overton $^3$}

\begin{document}

\maketitle

\footnotetext[1]{LAAS-CNRS, 7 Avenue du Colonel Roche, 31077
Toulouse, France}

\footnotetext[2]{Department of Control Engineering, Faculty of
Electrical Engineering, Czech Technical University in Prague,
Technick\'a 2, 166\,27 Prague, Czech Republic.  Email:
henrion@laas.fr. Research supported in part by Projects 102/06/0652
and 102/05/0011 of the Grant Agency of the Czech Republic and Project
ME 698/2003 of the Ministry of Education of the Czech Republic.}

\footnotetext[3]{Courant Institute of Mathematical Sciences,
New York University, 251 Mercer St., New York, NY, USA.
Email: overton@cs.nyu.edu. Research supported in part by 
National Science Foundation Grant DMS-0412049.}

\begin{abstract}
We consider the following problem: find a fixed-order linear 
controller that maximizes the closed-loop asymptotic decay rate for the 
classical two-mass-spring system.
This can be formulated as the problem of minimizing the
abscissa (maximum of the real parts of the roots)
of a polynomial whose coefficients depend linearly
on the controller parameters. We show that the only order for which 
there is a non-trivial solution is 2.  In this case, we 
derive a controller that we prove locally maximizes the asymptotic
decay rate,  
using recently developed techniques from nonsmooth analysis.
\end{abstract}

\section{Problem Statement}

\begin{figure}
\begin{center}
\includegraphics[scale=0.5]{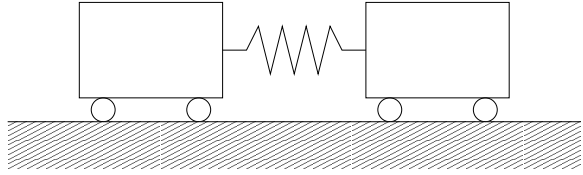}
\caption{Two-mass-spring system.\label{tms}}
\end{center}
\end{figure}

We consider the system shown in Figure \ref{tms}
consisting of two masses interconnected
by a spring, a typical control benchmark problem which is a generic
model of a system with a rigid body mode and one vibration mode
\cite{aiaa}.
If the first mass is pulled sufficiently far apart from the second mass and
suddenly dropped, then the two masses will oscillate until they reach
their equilibrium position.

The control problem we study in this note consists of 
appropriately moving the second mass so that the first mass settles down to
its final position as fast as possible; more specifically, we want to maximize
the asymptotic decay rate. For this we use a linear feedback
controller between the system output (measured position of the second
mass) and the system input (actuator positioning the first mass).
This control problem can be formulated as the minimization
of the abscissa (maximum of the real parts of the roots) of a polynomial 
in the complex Laplace
indeterminate, whose coefficients depend affinely on the controller
parameters. This polynomial is the denominator of the closed-loop system transfer
function.

For notational simplicity, we assume here that both
mass weights and the spring constant are normalized to one. As shown in
\cite{aiaa}, the polynomial is then given by
\[
p(x,y)(s) = (s^4+2s^2)x(s)+y(s)
\]
where $x(s)$ and $y(s)$ are respectively the denominator and numerator
polynomial of the controller transfer function to be determined.
This transfer function is assumed to be proper, i.e.,
\[
m = \mathrm{deg}\:x(s)\geq \mathrm{deg}\:y(s).
\]
The integer $m$ is called the order of the controller.
Without loss of generality, we take $x(s)$ to be monic.
Letting $\P^n$ denote the linear space of polynomials with complex\footnote{
We work with the space of polynomials with complex coefficients for technical reasons;
the chain rule that we use in Section \ref{locopt}
is most naturally stated in this context.}
coefficients and degree $\leq n$, we therefore write 
$y(s) \in \P^m$ and $x(s) - s^m \in \P^{m-1}$.

For a given polynomial $q \in \P^n$, we define the abscissa of $q$ by
\[
    \alpha(q) = \max \left \{ \Re~z: ~ q(z) = 0 \right \}.
\]
Since the roots of $p(x,y)$ are the closed-loop system poles,
the closed-loop two-mass-spring system is stable if and only if $\alpha(p(x,y)) <
0$. We are interested in maximizing the asymptotic decay rate of the system, 
i.e., solving the optimization problem
\begin{equation}\label{minabs}
\inf_{y(s) \in \P^m,~x(s) - s^m \in \P^{m-1}} \alpha(p(x,y)).
\end{equation}

Section \ref{first} shows that for order $m\leq 1$,
there is no stabilizing controller, i.e., no
 $(x,y)$ such that $\alpha(p(x,y)) < 0$.
Section \ref{third} shows that when $m \geq 3$,
the abscissa $\alpha(p(x,y))$ is unbounded below.
Section \ref{second} studies the more interesting case
$m=2$, and gives a formula for $(x,y)$ that, in Section \ref{locopt},
 we prove is a strong local minimizer of $\alpha(p(x,y))$.  
Section \ref{timefrag} plots the time response of the
optimized controller and discusses the issue of robustness.
Concluding remarks are made in Section \ref{conclusions}.

\section{First-Order Controller Design}\label{first}

If we assume that the controller has order $m=1$, then both $x(s)$ and
$y(s)$ are first degree polynomials, say
\[
x(s) = x_0 + s, \quad y(s) = y_0 + y_1 s,
\]
with
\[
p(x,y)(s) = y_0 + y_1 s + 2x_0 s^2 + 2 s^3 + x_0 s^4 + s^5.
\]
This polynomial is stable, i.e., it has all its roots in the open left
half-plane, if and only if all the principal minors of its Hurwitz
matrix
\[
\left[\begin{array}{ccccc}
x_0 & 1 & 0 & 0 & 0 \\
2x_0 & 2 & x_0 & 1 & 0 \\
y_0 & y_1 & 2x_0 & 2 & x_0 \\
0 & 0 & y_0 & y_1 & 2x_0 \\
0 & 0 & 0 & 0 & y_0
\end{array}\right]
\]
are all strictly positive, see e.g. \cite{dorato}.
This can never be the case since the 2-by-2 northwest
minor has rank one for all $x_0$. Hence a controller of first order
(or less) cannot stabilize the two-mass-spring system.

Note that in problem (\ref{minabs}) the minimum abscissa
$\alpha(p(x,y)) = 0$ is attained for any static controller $x(s)=1$,
$y(s)=y_0=k$ with $k \in [0,1]$, because then
\[
p(x,y)(s) = k+2s^2+s^4 = (s^2+1+\sqrt{1-k})(s^2+1-\sqrt{1+k})
\]
has only imaginary roots.
Similarly, $\alpha(p(x,y)) = 0$ is attained for any first-order
controller $x(s) = x_0 + s$, $y(s) = y_0 + y_1 s$ such that $x_0 =
y_0 = 0$, $y_1 = k$ and $k \in [0,1]$, since then
\[
p(x,y)(s) = (k+2s^2+s^4)s.
\]

\section{Third-Order Controller Design}\label{third}

If we seek a controller of order $m=3$ then we can write
\[
p(x,y)(s) = (s^4+2s^2)(x_0+x_1s+x_2s^2+s^3)+y_0+y_1s+y_2s^2+y_3s^3
= s^7 + \sum_{i=0}^6 p_i s^i.
\]
By identifying powers of the
indeterminate $s$, we derive the linear system of equation
\[
\left[\begin{array}{ccccccc}
0 & 0 & 0 & 1 & 0 & 0 & 0 \\
0 & 0 & 0 & 0 & 1 & 0 & 0 \\
2 & 0 & 0 & 0 & 0 & 1 & 0 \\
0 & 2 & 0 & 0 & 0 & 0 & 1 \\
1 & 0 & 2 & 0 & 0 & 0 & 0 \\
0 & 1 & 0 & 0 & 0 & 0 & 0 \\
0 & 0 & 1 & 0 & 0 & 0 & 0
\end{array}\right] 
\left[\begin{array}{c}x_0\\x_1\\x_2\\y_0\\y_1\\y_2\\y_3\end{array}\right]
= \left[\begin{array}{c}p_0\\p_1\\p_2\\p_3\\p_4\\p_5-2\\p_6 \end{array}\right].
\]
This 7-by-7 matrix is called the Sylvester matrix of polynomials
$s^4+2s^2$ and $1$, and it is non-singular because these two
polynomials share no common roots. In other words, we can find
controller coefficients defining any desired closed-loop 
polynomial. This is what Dorato \cite{dorato} calls the fundamental
theorem of feedback control, namely the fact that the
poles of a single-input-single-output linear system of order $n$ can
be placed arbitrarily by a linear controller of order $n-1$.
Consequently, an arbitrarily large negative abscissa can be achieved
in closed-loop by a controller of order three (or more).

For example, by solving the linear system of equations shown above, we
obtain that the controller polynomials
\[
\begin{array}{rcl}
x(s) & = & (-35z^3+14z)+(21z^2-2)s-7z s^2+s^3, \\
y(s) & = & -z^7+7z^6 s+(-21z^5+70z^3-28z)s^2
+(35z^4-42z^2+4)s^3
\end{array}
\]
place all the closed-loop poles at an arbitrary real value
$z$.  It follows that problem (\ref{minabs}) is not bounded below for
$m \geq 3$.

\section{Second-Order Controller Design}\label{second}

In the case of a controller of order $m=2$, we have
\begin{equation}
\label{pxydef}
p(x,y)(s) = (s^4+2s^2)(x_0+x_1s+s^2)+y_0+y_1s+y_2s^2.
\end{equation}
We can cluster all the closed-loop poles at a real
negative value $z$ by solving the following system
of equations
\[
\begin{array}{rcl}
y_0 & = & z^6 \\
y_1 & = & -6z^5 \\
y_2 + 2x_0 & = & 15z^4 \\
2x_1 & = & -20z^3 \\
2+x_0 & = & 15z^2 \\
x_1 & = & -6z.
\end{array}
\]
We observe that the only constraint on $z$ is enforced
by the fourth and sixth equations, namely
\[
5z^3 = 3z.
\]
We rule out the case $z=0$ since the closed-loop system
would be only marginally stable, and we extract the
negative solution
\[
z^* = -\frac{\sqrt{15}}{5} \approx -0.7746.
\]
The controller coefficients can now be derived by substitution, resulting in
\begin{equation}
\label{xystar}
x_0^*=7,~ x_1^* = \frac{6\sqrt{15}}{5},~y_0^* = \frac{27}{125},
  ~y_1^* = \frac{54\sqrt{15}}{125},~y_2^*= -\frac{43}{5},
\end{equation}
yielding 
\begin{equation}\label{pxystar}
p(x^*,y^*)(s) = s^6+\frac{6\sqrt{15}}{5}s^5+9s^4+\frac{12\sqrt{15}}{5}s^3+\frac{27}{5}s^2
+\frac{54\sqrt{15}}{125}s+\frac{27}{125}.
\end{equation}
We should note that we realized that all roots could be clustered at
a single point $z^*$ after performing numerical experiments using {\sc hifoo}
\cite[Section 6]{rocond06}, a new toolbox for low order controller design using
methods of nonsmooth optimization.

\section{Local Optimality Certificate}\label{locopt}

In this section we prove that $z^* = -\frac{\sqrt{15}}{5}$ is,
at least locally, the minimal abscissa that can be achieved with
a second-order controller, i.e., a local minimizer of problem
(\ref{minabs}) for $m=2$.  
This is nontrivial, since one might think it is necessary to consider
all possible splittings of the multiple root under perturbation.
We prove local optimality using recent advances in nonsmooth analysis.

Recall that $\P^n$ is the linear space of complex polynomials with degree less than 
or equal to $n$ (with complex dimension $n+1$) and let $\C^n$ denote the space of complex
vectors of length $n$. We write elements of $\C^n$ as row vectors.

We show that $(x^*,y^*)$ defined by (\ref{xystar})
locally optimizes the abscissa of $p(x,y)$, in the sense that 
any sufficiently small perturbation to (\ref{xystar}) strictly increases the 
maximum of the real parts of the roots.  In fact, we prove that (\ref{xystar}) 
is a sharp local minimizer, in the following sense.  

\medskip

{\bf Theorem} 
\emph{The abscissa of the polynomial $p(x,y)$ defined in (\ref{pxydef})
is locally minimized by $(x^*,y^*)$ with coefficients given by (\ref{xystar}).  
Furthermore, for $(x,y)$ sufficiently close to $(x^*,y^*)$, we have
\[
         \alpha(p(x,y)) \geq \alpha(p(x^*,y^*)) + \tau \|d\|
\]
where $\tau$ is a positive constant and
\[ 
  d = \left [ x_0 - x_0^*,  x_1 - x_1^*, y_0 - y_0^*, y_1 - y_1^*, y_2 - y_2^* \right ].
\]
}

The proof of this theorem is the subject of the rest of this section.
It follows \cite[Section III A]{BurHenLewOveChoc} quite closely. 
We start by making a change of variables 
to the polynomial $p(x,y)(s)$, namely
{\samepage
\begin{eqnarray*}
               t &=& s - z^*,  \\
               q &=& [q_0,~q_1] = [x_0,~ x_1] - [x_0^*,~x_1^*], \\
               r &=& [r_0,~r_1,~r_2] = [y_0,~y_1,~y_2] - [y_0^*,~y_1^*~~y_2^*].
\end{eqnarray*}
}
A few lines of {\sc Maple} \cite{maple}, specifically
\begin{verbatim}
p:=(s^4+2*s^2)*(s^2+x1*s+x0)+(y2*s^2+y1*s+y0);
subs(x0=q0+7,x1=q1+6*sqrt(15)/5,y0=r0+27/125, y1=r1+54*sqrt(15)/125, 
   y2=r2-43/5,s=t+a,a=-sqrt(15)/5,p);
collect(%,t);
simplify(%);
collect(%,t)
\end{verbatim}
show that the resulting polynomial is 
\begin{equation}
        t \mapsto t^6 + A(q,r)(t)
\label{newpoly}
\end{equation}
where the linear map $A:\C^5\rightarrow \P^5$ is given by
{\samepage
\begin{eqnarray*}
A(q,r)(t) &=&
 q_1 t^5+\left (q_0-\sqrt{15}\, q_1\right) t^4+\left (8 q_1 -\frac{4}{5} \sqrt{15}\, q_0\right ) t^3+ \\
&&  \left (\frac{28}{5} q_0 -\frac{12}{5} \sqrt{15}\, q_1 + r_2\right ) t^2 + \\
&&  \left (\frac{-32}{25} \sqrt{15}\, q_0 +\frac{27}{5} q_1+r_1-\frac{2}{5} \sqrt{15}\, r_2\right ) t + \\
&& \frac{39}{25} q_0 -\frac{39}{125} \sqrt{15}\, q_1 +r_0 -\frac{1}{5} \sqrt{15}\, r_1 +\frac{3}{5} r_2.
\end{eqnarray*}
}
It is easily verified that $A(0,0)=0$; hence the map $A$ is indeed linear.  
Clearly, minimizing the abscissa of the polynomial (\ref{newpoly})
over $[q,~r] \in \C^5$ is equivalent to the original problem (\ref{minabs}).
Because the space of monic polynomials is not a linear space, it is convenient to
introduce the notation
\[
      \gamma(w) = \max\{\Re~t: t^{n+1} + w(t) = 0\}, \quad w \in \P^n,
\]
for the abscissa of  $t^{n+1} + w(t)$.  
We wish to establish that $0$ is a sharp local minimizer
of the composition of the function $\gamma$ with the linear map
$A$ over $[q,~r]$ in the parameter space $\C^5$.

To proceed further we need the notion of the adjoint map
$A^*:  \P^5 \rightarrow \C^5$, defined by 
\[
     \Big \langle w(t),~ A(q,r)(t)  \Big \rangle = \Big \langle A^*(w), ~[q,~r] \Big \rangle,
\]
for all polynomials $w \in \P^5$ and vectors $[q,~r] \in \C^5$, where the second inner product
is the usual real inner product on $\C^5$ and the first is a real inner product on
$\P^5$, namely 
\[
        \Big \langle \sum_{j=0}^5 c_j t^j, ~\sum_{j=0}^5 d_j t^j \Big \rangle = 
          \Re \sum_{j=0}^5 c_j \bar{d}_j.
\]
It is easy to see that $A^*$ is given by
\begin{equation}
\label{Astar}
  A^*\left( \sum_{j=0}^5 c_j t^j \right) = \left [ \begin{array}{c}
           \frac{39}{25} c_0 - \frac{32}{25}\sqrt{15}\, c_1 + \frac{28}{5} c_2 
          -\frac{4}{5}\sqrt{15}\, c_3 + c_4 \\
           -\frac{39}{125}\sqrt{15}\, c_0 + \frac{27}{5} c_1 - \frac{12}{5}\sqrt{15}\, c_2
              +8 c_3 -\sqrt{15}\, c_4 + c_5 \\
            c_0 \\
            - \frac{1}{5}\sqrt{15}\, c_0 + c_1 \\
            \frac{3}{5} c_0 - \frac{2}{5}\sqrt{15}\, c_1 + c_2
             \end{array} \right ] .
\end{equation}
Following \cite{BurHenLewOveChoc} and \cite{BurLewOveAMS},
we will establish that 0 is a sharp 
local minimizer of the composition of $\gamma$ with the linear map $A$, which we denote
$\gamma \circ A$, by showing that
\begin{equation}
\label{intpartial}
      0 \in \interior ~\partial \left ( \gamma \circ A \right ) (0),
\end{equation}
where $\partial$ is the subdifferential operator of variational analysis
 \cite{BurLewOveAMS}, \cite[Chap.\ 8]{RocWet98}. In order to do this we can use 
the nonsmooth chain rule \cite[Lemma 4.4]{BurLewOveAMS}
\begin{equation}
\label{chainrule}
       \partial \left ( \gamma \circ A \right ) = A^* \partial \gamma(0), 
\end{equation}
as long as we verify the constraint qualification
\begin{equation}
\label{conqual}
      \N(A^*) \cap \partial^\infty \gamma(0) = \{ 0 \},
\end{equation}
where $\N$ denotes null space and $\partial^\infty$ denotes the horizon 
subdifferential operator  \cite{BurLewOveAMS}, \cite[Chap.\ 8]{RocWet98}.
This chain rule is valid because of the \emph{subdifferential regularity}
\cite[Chap.\ 8]{RocWet98}
of the function $\gamma$ on $\P^n$, established in \cite{BurOvePoly}.
The following formulas for the subdifferential and horizon subdifferential
of $\gamma$ at 0 were also established in \cite{BurOvePoly}, but we 
follow the notation used in \cite[Theorem 3.3]{BurHenLewOveChoc}: 
\begin{eqnarray*}
     \partial \gamma(0) &=&
        \Big\{ \sum_{j=0}^n c_j t^j : c_n = -\frac{1}{n+1},~\Re~c_{n-1} \le 0 \Big\},
 \\ 
 \partial^{\infty} \gamma(0) &=&
    \Big\{ \sum_{j=0}^n c_j t^j : c_n = 0,~\Re~c_{n-1} \le 0 \Big\}. 
\end{eqnarray*}
It follows from the latter formula that the constraint qualification (\ref{conqual}) 
holds if $c_5=0$ and $ A^*\left (\sum_{j=0}^5 c_j t^j \right ) = 0$ implies 
$c=[c_0,\cdots,c_5] = 0$, a fact that
is easily checked by observing that the 2 by 2 linear system
\[
     \left [ \begin{array}{cc} 
                  -\frac{4}{5}\sqrt{15} & 1 \\
                  8 & -\sqrt{15}  \end{array} \right ]
     \left [ \begin{array}{c} c_3 \\ c_4 \end{array} \right ]
     = \left [ \begin{array}{c} 0 \\ 0 \end{array} \right ]
\]
has only the trivial solution $[c_3,~c_4] = 0$.
Thus the nonsmooth chain rule (\ref{chainrule}) yields
\[
       \partial \left ( \gamma \circ A \right ) = 
     \{ A^* \left ( \sum_{j=0}^5 c_j t^j \right) : c_5 = -\frac{1}{6}, ~~ \Re~c_4 \leq 0 \}.
\]
The final step is to determine whether $0$ is in the interior of this subdifferential set. 
To check this, we need to solve the following linear system: set the right-hand side of 
(\ref{Astar}) to 0 as well as $ c_5 = -\frac{1}{6}$, which reduces to
\[
     \left [ \begin{array}{cc} 
                  -\frac{4}{5}\sqrt{15} & 1 \\
                  8 & -\sqrt{15}  \end{array} \right ]
     \left [ \begin{array}{c} c_3 \\ c_4 \end{array} \right ]
     = \left [ \begin{array}{c} 0 \\ \frac{1}{6} \end{array} \right ].
\]
This linear system has a unique solution, namely 
$[c_3 ,~ c_4] = [-\frac{1}{24}~~-\frac{1}{30}\sqrt{15}\,]$.  Since this satisfies the
inequality $\Re~c_4 \leq 0$, it follows that $0$ is in the subdifferential set, and 
furthermore, since the inequality holds strictly,
that every point near $0$ is in the subdifferential set, and therefore
that (\ref{intpartial}) holds.  This completes the proof of the theorem.

\section{Time Response and Fragility}\label{timefrag}

\begin{figure}
\begin{center}
\includegraphics[scale=0.7]{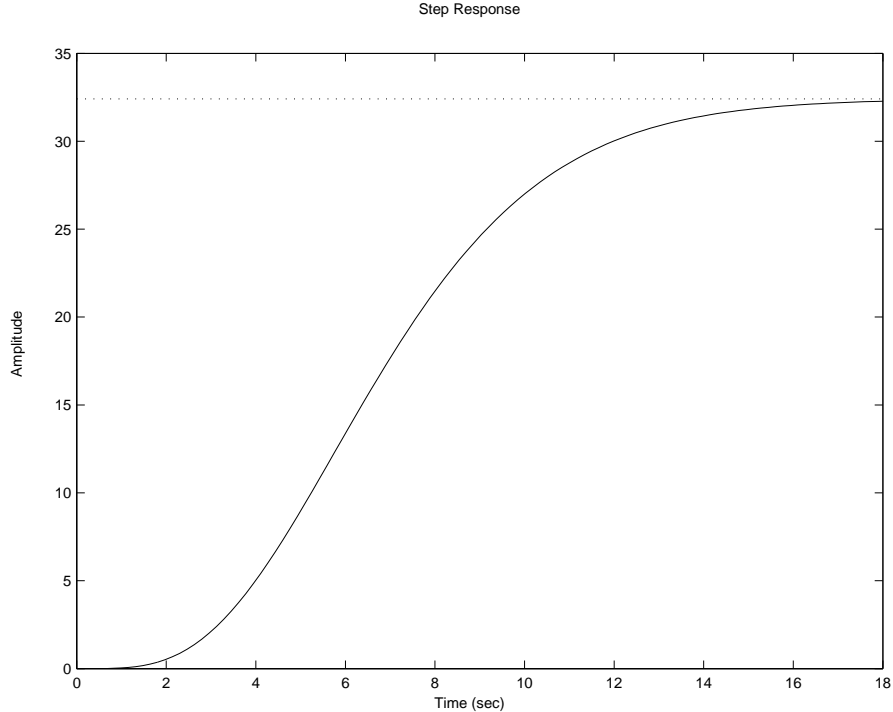}
\caption{Closed-loop step response.\label{step6}}
\end{center}
\end{figure}

In Figure \ref{step6} we graph the step response of the two mass-spring
system fed back with controller (\ref{xystar}), obtained with the
following {\sc Matlab} \cite{matlab} commands:
\begin{verbatim}
P = tf(1,[1 0 2 0 0]);
K = tf([-43/5 54*sqrt(15)/125 27/125],[1 6*sqrt(15)/5 7]);
T = feedback(P,K);
step(T);
\end{verbatim}
We can see that the settling time is around 16 seconds.

It is well known that multiple roots of polynomials are very sensitive to
perturbations in the coefficients. In practice, this means that the
closed-loop system will be fragile, or non-robust, with respect to
uncertain data, implementation errors, or even rounding errors. For
example, if instead of implementing
the exact second-order controller (\ref{xystar}),
we implement the nearby controller, obtained by keeping 5
significant digits, given by
\[
x_0=7,~ x_1 = 4.6476,~y_0 = 0.2160,~y_1 =1.6731,~y_2=-8.6,
\]
then we obtain closed-loop poles at $-0.9405$, $-0.8163\pm 0.1489~i$,
$-0.7500$ and $-0.6622 \pm 0.0786~i$, quite far from the single
pole at $-0.7746$ assigned with the exact controller.

\begin{figure}
\begin{center}
\includegraphics[scale=0.7]{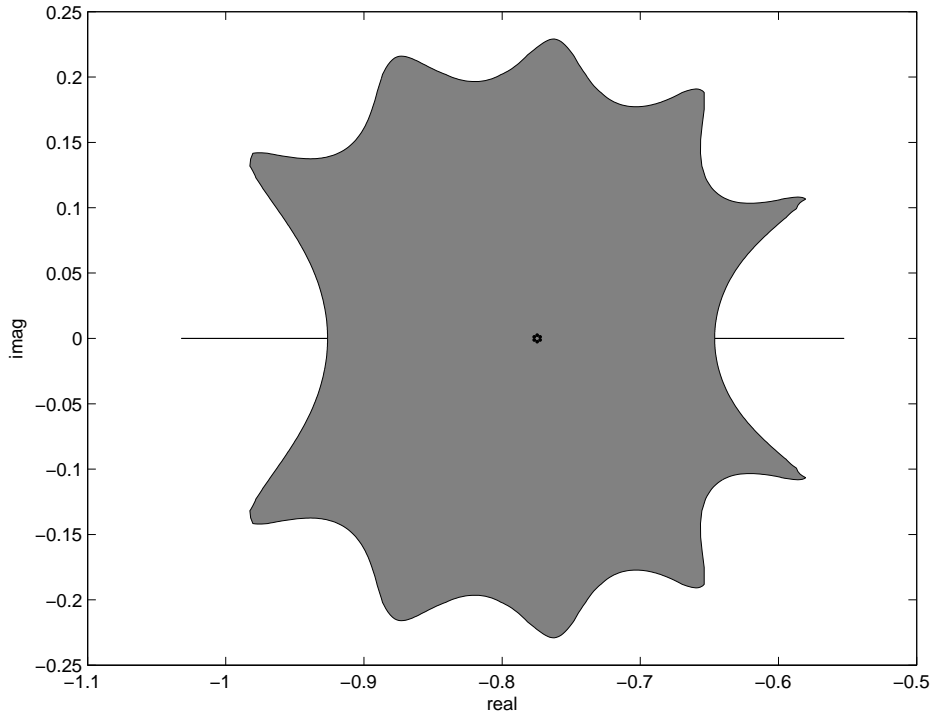}
\caption{Possible zeros of the closed-loop
polynomial $p(x^*,y^*)$ under real perturbations to its coefficients
with norm $\leq 10^{-4}$.\label{rps6}}
\end{center}
\end{figure}

This phemonenon can be studied graphically.
In Figure \ref{rps6} we show\footnote{Thanks to S.\ Graillat,
  N.\ Higham and F.\ Tisseur for providing {\sc matlab} scripts
  for the computation of real pseudozero sets of a polynomial.}, in the gray
region, all possible roots of polynomials that can be obtained by
real perturbations to the polynomial $p(x^*,y^*)$ given in (\ref{pxystar}),
where the norm of the vector of perturbations to the coefficients is
no more than $\epsilon=10^{-4}$.  This region is sometimes called
the real pseudozero set.

\section{Concluding Remarks}\label{conclusions}

In this note we have formulated the problem of maximizing the
closed-loop asymptotic decay rate of a linear control system
as a nonsmooth, nonconvex problem of polynomial abscissa minimization,
focusing on the case of a benchmark two-mass-spring system. We derived
a formula for a second-order controller with closed-loop poles 
clustered at a single point.  Our main contribution
is the use of recently developed techniques from nonsmooth variational
analysis to prove local optimality of this controller.

Motivated by this result, as well as the result in \cite{BurHenLewOveChoc}
on which it is based, very recent work \cite{BloOve06} using a completely
different technique shows that the second-order controller described
above is actually globally optimal.

Finally, it should be emphasized that asymptotic decay rate
maximization is not, by itself, a practical objective.
As shown graphically in Figure \ref{rps6}, our
locally optimal controller yields a closed-loop system which is
sensitive to uncertainty and/or disturbance. In other
words, lack of robustness is the price one has to pay to maximize the
decay rate. In a typical control engineering system, a trade-off
should be found between the asymptotic decay rate and other quantities, such
as the complex or real stability radius, the complex or 
real pseudoabscissa (the maximum real part of the points in the complex 
or real pseudozero set), or $H_2$ or $H_{\infty}$ performance measures.

%

\end{document}